\documentclass[12pt]{amsart}
\usepackage{amssymb}
\usepackage[all]{xy}

\begin{document}
\theoremstyle{plain}
\newtheorem{thm}{Theorem}[section]
\newtheorem*{thm1}{Theorem 1}
\newtheorem*{thm1.1}{Theorem 1.1}
\newtheorem*{thmM}{Main Theorem}
\newtheorem*{thmA}{Theorem A}
\newtheorem*{thm2}{Theorem 2}
\newtheorem{lemma}[thm]{Lemma}
\newtheorem{lem}[thm]{Lemma}
\newtheorem{cor}[thm]{Corollary}
\newtheorem{pro}[thm]{Proposition}
\newtheorem{propose}[thm]{Proposition}
\newtheorem{variant}[thm]{Variant}
\theoremstyle{definition}
\newtheorem{notations}[thm]{Notations}
\newtheorem{rem}[thm]{Remark}
\newtheorem{rmk}[thm]{Remark}
\newtheorem{rmks}[thm]{Remarks}
\newtheorem{defi}[thm]{Definition}
\newtheorem{exe}[thm]{Example}
\newtheorem{claim}[thm]{Claim}
\newtheorem{ass}[thm]{Assumption}
\newtheorem{prodefi}[thm]{Proposition-Definition}
\newtheorem{que}[thm]{Question}
\newtheorem{con}[thm]{Conjecture}

\newtheorem*{dmlcon}{Dynamical Mordell-Lang Conjecture}
\newtheorem*{condml}{Dynamical Mordell-Lang Conjecture}
\numberwithin{equation}{section}
\newcounter{elno}                % This to number lists
\def\points{\list
{\hss\llap{\upshape{(\roman{elno})}}}{\usecounter{elno}}}
\let\endpoints=\endlist

\newcommand{\la}{\lambda}
\newcommand{\mc}{\mathcal}
\newcommand{\mb}{\mathbb}
\newcommand{\surj}{\twoheadrightarrow}
\newcommand{\inj}{\hookrightarrow}
\newcommand{\zar}{{\rm zar}}
\newcommand{\an}{{\rm an}}
\newcommand{\red}{{\rm red}}
\newcommand{\codim}{{\rm codim}}
\newcommand{\Supp}{{\rm Supp}}
\newcommand{\rank}{{\rm rank}}
\newcommand{\Ker}{{\rm Ker \ }}
\newcommand{\Pic}{{\rm Pic}}
\newcommand{\Div}{{\rm Div}}
\newcommand{\Hom}{{\rm Hom}}
\newcommand{\im}{{\rm im}}
\newcommand{\Spec}{{\rm Spec \,}}
\newcommand{\Nef}{{\rm Nef \,}}
\newcommand{\Frac}{{\rm Frac \,}}
\newcommand{\Sing}{{\rm Sing}}
\newcommand{\sing}{{\rm sing}}
\newcommand{\reg}{{\rm reg}}
\newcommand{\Char}{{\rm char}}
\newcommand{\Tr}{{\rm Tr}}
\newcommand{\ord}{{\rm ord}}
\newcommand{\id}{{\rm id}}
\newcommand{\NE}{{\rm NE}}
\newcommand{\Gal}{{\rm Gal}}
\newcommand{\Min}{{\rm Min \ }}
\newcommand{\Max}{{\rm Max \ }}
\newcommand{\Alb}{{\rm Alb}\,}
\newcommand{\GL}{{\rm GL}\,}        % For the general linear group
\newcommand{\PGL}{{\rm PGL}\,}
\newcommand{\Bir}{{\rm Bir}}
\newcommand{\Aut}{{\rm Aut}}
\newcommand{\End}{{\rm End}}
\newcommand{\Per}{{\rm Per}\,}
\newcommand{\ie}{{\it i.e.\/},\ }
\newcommand{\niso}{\not\cong}
\newcommand{\nin}{\not\in}
\newcommand{\soplus}[1]{\stackrel{#1}{\oplus}}
\newcommand{\by}[1]{\stackrel{#1}{\rightarrow}}
\newcommand{\longby}[1]{\stackrel{#1}{\longrightarrow}}
\newcommand{\vlongby}[1]{\stackrel{#1}{\mbox{\large{$\longrightarrow$}}}}
\newcommand{\ldownarrow}{\mbox{\Large{\Large{$\downarrow$}}}}
\newcommand{\lsearrow}{\mbox{\Large{$\searrow$}}}
\renewcommand{\d}{\stackrel{\mbox{\scriptsize{$\bullet$}}}{}}
\newcommand{\dlog}{{\rm dlog}\,}    % For dlog
\newcommand{\longto}{\longrightarrow}
\newcommand{\vlongto}{\mbox{{\Large{$\longto$}}}}
\newcommand{\limdir}[1]{{\displaystyle{\mathop{\rm lim}_{\buildrel\longrightarrow\over{#1}}}}\,}
\newcommand{\liminv}[1]{{\displaystyle{\mathop{\rm lim}_{\buildrel\longleftarrow\over{#1}}}}\,}
\newcommand{\norm}[1]{\mbox{$\parallel{#1}\parallel$}}
\newcommand{\boxtensor}{{\Box\kern-9.03pt\raise1.42pt\hbox{$\times$}}}
\newcommand{\into}{\hookrightarrow}
\newcommand{\image}{{\rm image}\,}
\newcommand{\Lie}{{\rm Lie}\,}      % For Lie algebra of groups
\newcommand{\CM}{\rm CM}
\newcommand{\sext}{\mbox{${\mathcal E}xt\,$}}  % For sheaf Ext
\newcommand{\shom}{\mbox{${\mathcal H}om\,$}}  %For sheaf Hom
\newcommand{\coker}{{\rm coker}\,}  % For the cokernel of a morphism
\newcommand{\sm}{{\rm sm}}
\newcommand{\pgcd}{\text{pgcd}}
\newcommand{\trd}{\text{tr.d.}}
\newcommand{\tensor}{\otimes}
\renewcommand{\iff}{\mbox{ $\Longleftrightarrow$ }}
\newcommand{\supp}{{\rm supp}\,}
\newcommand{\ext}[1]{\stackrel{#1}{\wedge}}
\newcommand{\onto}{\mbox{$\,\>>>\hspace{-.5cm}\to\hspace{.15cm}$}}
\newcommand{\propsubset}
{\mbox{$\textstyle{
\subseteq_{\kern-5pt\raise-1pt\hbox{\mbox{\tiny{$/$}}}}}$}}
% Skriptbuchstaben
\newcommand{\sA}{{\mathcal A}}
\newcommand{\sB}{{\mathcal B}}
\newcommand{\sC}{{\mathcal C}}
\newcommand{\sD}{{\mathcal D}}
\newcommand{\sE}{{\mathcal E}}
\newcommand{\sF}{{\mathcal F}}
\newcommand{\sG}{{\mathcal G}}
\newcommand{\sH}{{\mathcal H}}
\newcommand{\sI}{{\mathcal I}}
\newcommand{\sJ}{{\mathcal J}}
\newcommand{\sK}{{\mathcal K}}
\newcommand{\sL}{{\mathcal L}}
\newcommand{\sM}{{\mathcal M}}
\newcommand{\sN}{{\mathcal N}}
\newcommand{\sO}{{\mathcal O}}
\newcommand{\sP}{{\mathcal P}}
\newcommand{\sQ}{{\mathcal Q}}
\newcommand{\sR}{{\mathcal R}}
\newcommand{\sS}{{\mathcal S}}
\newcommand{\sT}{{\mathcal T}}
\newcommand{\sU}{{\mathcal U}}
\newcommand{\sV}{{\mathcal V}}
\newcommand{\sW}{{\mathcal W}}
\newcommand{\sX}{{\mathcal X}}
\newcommand{\sY}{{\mathcal Y}}
\newcommand{\sZ}{{\mathcal Z}}
% Sonderbuchstaben mit Doppellinie
\newcommand{\A}{{\mathbb A}}
\newcommand{\B}{{\mathbb B}}
\newcommand{\C}{{\mathbb C}}
\newcommand{\D}{{\mathbb D}}
\newcommand{\E}{{\mathbb E}}
\newcommand{\F}{{\mathbb F}}
\newcommand{\G}{{\mathbb G}}
\newcommand{\HH}{{\mathbb H}}
\newcommand{\I}{{\mathbb I}}
\newcommand{\J}{{\mathbb J}}
\newcommand{\M}{{\mathbb M}}
\newcommand{\N}{{\mathbb N}}
\renewcommand{\P}{{\mathbb P}}
\newcommand{\Q}{{\mathbb Q}}
\newcommand{\R}{{\mathbb R}}
\newcommand{\T}{{\mathbb T}}
\newcommand{\U}{{\mathbb U}}
\newcommand{\V}{{\mathbb V}}
\newcommand{\W}{{\mathbb W}}
\newcommand{\X}{{\mathbb X}}
\newcommand{\Y}{{\mathbb Y}}
\newcommand{\Z}{{\mathbb Z}}

\newcommand{\fix}{\mathrm{Fix}}

%%%%%%%%%%%%%%%%%%%%%%%%%%%%%%%%%%%%%%%%%%%%%%%%%%%%%%%%%%%%%%
%%%%%%%%%%%%%%%%%%%%%%%%%%%%%%%%%%%%%%%%%%%%%%%%%%%%%%%%%%%%%%
\title[]{Dynamical Mordell-Lang conjecture for birational polynomial morphisms on $\mathbb{A}^2$}

\author{Xie Junyi}
%%\address{\'Ecole Normale Sup\'erieure, 45 rue d'Ulm, 75005, Paris, France
%%}
\address{Centre de Math\'ematiques Laurent Schwartz \'Ecole Polytechnique, 91128, Palaiseau
Cedex, France}

\email{junyi.xie@math.polytechnique.fr}
\date{\today}
\bibliographystyle{plain}
\maketitle

\begin{abstract}
We prove the dynamical Mordell-Lang conjecture for birational polynomial morphisms on $\mathbb{A}^2$.
\end{abstract}
\renewcommand{\thefootnote}{}
\footnotetext{The author is supported
by the ERC-starting grant project ¡±Nonarcomp¡± no.307856.}

\section{Introduction}
The Mordell-Lang conjecture proved by Faltings \cite{Faltings1994} and Vojta \cite{Vojta1996} says that if $V$ is
a subvariety of a semiabelian variety $G$ defined over $\mathbb{C}$ and $\Gamma$ is a finitely generated subgroup  of $G(\mathbb{C})$, then $V(\mathbb{C})\bigcap \Gamma$ is a union of at most finitely many translates of
subgroups of $\Gamma$.
\smallskip

The following dynamical analogue of the Mordell-Lang conjecture was proposed by
Ghioca and Tucker.
\begin{dmlcon}[\cite{Ghioca2009}]\label{dml}Let $X$ be a quasiprojective variety defined over $\mathbb{C}$, let
$f: X \rightarrow X$ be an endomorphism, and $V$ be any subvariety of $X$. For any point $p\in X(\mathbb{C})$  the set $\{n\in \mathbb{N}|\,\,f^n(p)\in V(\mathbb{C})\}$ is a union of at most finitely many arithmetic progressions.
\end{dmlcon}
An arithmetic progression is a set of the form $\{an+b|\,\, n\in \mathbb{N}\}$ with $a,b\in \mathbb{N}$ possibly with $a=0$.

Observe that this conjecture implies the classical Mordell-Lang conjecture in the case $\Gamma\simeq (\mathbb{Z},+)$.

\smallskip

The Dynamical Mordell-Lang conjecture has been proved by Denis \cite{Denis1994} for automorphisms of projective spaces and was later generalized by Bell \cite{Bell2006} to the case of automorphisms of affine varieties. In \cite{Bell2010}, Bell, Ghioca and Tucker proved it for \'etale maps of quasiprojective varieties. The conjecture is also known in the case where $f=( F(x_1) , G(x_2)):\mathbb{A}^2_{\mathbb{C}}\rightarrow \mathbb{A}^2_{\mathbb{C}}$ where
$F,G$ are polynomials and the subvariety $V$ is a line (\cite{Ghioca2008}), and in the case
$f=(F(x_1),\cdots, F(x_n)):\mathbb{A}^n_{K}\rightarrow \mathbb{A}^n_{K}$
where $F\in K[t]$ is an indecomposable polynomial defined over a number field $K$ which
has no periodic critical points other than the point at infinity and $V$ is a curve (\cite{Benedetto2012}).

\bigskip

Our main result can be stated as follows.

\begin{thmA}\label{mainpolyarb}
Let $K$ be any algebraically closed field of characteristic $0$, and $f:\mathbb{A}^2_K\rightarrow\mathbb{A}^2_K$ be any birational polynomial morphism defined over $K$. Let $C$ be any curve in $\mathbb{A}^2_K$, and $p$ be any point in $\mathbb{A}^2(K)$. Then the set $\{n\in \mathbb{N}|\,\, f^n(p)\in C\}$ is a union of at most finitely many arithmetic progressions.
\end{thmA}
In the case the map is an automorphism of $\mathbb{A}^2_K$ of H\'enon type (see
\cite{Friedland1989}) then this result follows from \cite{Bell2010}. Our proof provides however
an alternative approach and does not rely on the construction of $p$-adic invariant
curves. \medskip

Recall that the $algebraic$ $degree$  of a polynomial transformation $f(x,y)=(f_1(x,y),f_2(x,y))$ is defined by  $\deg f:= \max\{\deg f_1,\deg f_2\}$.
The limit  $\la (f):=\lim_{n\rightarrow \infty}(\deg f^n)^{1/n}$ exists and we refer to it as the $dynamical$ $degree$ of $f$ (see \cite{favre,Dinh2005}). Our proof shows that when $\la (f)>1$, then Theorem A holds for fields of arbitrary characteristic.

Note however that our Theorem A does not hold when $\Char K>0$ and $\lambda(f)=1$ (see \cite[Proposition 6.1]{Bell2006} for a counter-example).
\bigskip

To explain our strategy, we fix a birational polynomial morphism $f:\mathbb{A}_K^2\rightarrow \mathbb{A}_K^2$. By some reduction arguments, we may assume that $K=\overline{\mathbb{Q}}$.

We may compactify $\mathbb{A}^2$ by \cite{Favre2011} to a smooth projective surface, such that $f$ extends to a birational transformation on $X$ fixing a point $Q$ in $X\setminus \mathbb{A}^2$, and $f$ contracts all curves at infinity to $Q$ (see \cite{Favre2011} and Section \ref{l=1model}).

The key idea of our proof is to take advantage of this attracting fixed point and to apply the
following local version of the Dynamical Mordell-Lang conjecture.
\begin{thm}\label{local}Let $X$ be a smooth projective surface over an arbitrary valued field $(K,|\cdot|)$  and $f: X\dashrightarrow X$ be a birational transformation defined over $K$. Let $C$ be any curve in $X$. Pick any $K$-point $p$ such that $f^n(p)\in X\setminus I(f)$ for all integers $n\geq 0$, and $f^n (p)$ tends to a fixed $K$-point $Q\in I(f^{-1})\setminus I(f)$ with respect to a projective metric induced by $|\cdot|$ on $X$.

If the set $$\{n\in \mathbb{N}|\,\, f^n(p)\in C\}$$ is infinite, then either $f^n(p)=Q$ for some $n\geq 0$ or $C$ is fixed.
\end{thm}

To complete the proof of Theorem A we now rely on a global argument.
When the curve $C$ is passing through the fixed point $Q$ in $X$,
we cover the $\overline{\mathbb{Q}}$-points of the curve $C$ by the
basin of attraction of $Q$ with respect to all absolute values
on $\overline{\mathbb{Q}}$. If the point $p$ belongs to one of these attracting basins, then the local dynamical Mordell-Lang
applies and we are done. Otherwise it is possible
to bound the height of $p$ and Northcott theorem shows that it is periodic.

Finally when neither the curve $C$ nor its iterates contain the fixed point $Q$, we are in
position to apply the next result which allows us to conclude.
\begin{thm}\label{qinf^nC}Let $X$ be a smooth projective surface over an algebraically closed field, $f:X\dashrightarrow X$ be an algebraically stable birational transformation and $C$ be an irreducible curve in $X$ such that $f^n$ does not contract $C$ for any $n$.

If $f^n(C)\bigcap I(f)\neq \emptyset$ for all $n$, then $C$ is periodic.
\end{thm}
\bigskip

The article is organized in 8 sections.
In Section \ref{senoaba} we give background informations on birational surface maps and metrics on projective varieties defined over a valued field.
In Section \ref{acf} we prove Theorem \ref{qinf^nC}, which is a criterion for a curve to be periodic.
In Section \ref{sedml} we prove some basic properties for the maps satisfying the conclusion of dynamical Mordell-Lang conjecture.
In Section \ref{seldml} we prove Theorem \ref{local}.
In Section \ref{l=1} we prove Theorem A in the case the dynamical degree $\la(f)=1.$
In Section \ref{seupheight} we prove a technical lemma which gives a upper bound on height when $\la(f)>1$.
In Section \ref{seproofa} we prove Theorem A.

\section*{Acknowledgement}
I would like to thank
C. Favre for his support and his direction during the
writing of this article.

\section{Notations and basics}\label{senoaba}
\subsection{Basics on birational maps on surfaces}\label{bobmos}See \cite{debarre,favre,Favre2011} for details.

In this section a variety is defined over an algebraically closed field $k$.
Recall that the resolution of singularities exists for surfaces over any algebraically closed field (see \cite{Abhyankar1966}).

Let $X$ be a smooth projective surface. We denote by $N^1(X)$ the N\'eron-Severi group of $X$ i.e. the group of numerical
equivalence classes of divisors on $X$ and write $N^1(X)_{\mathbb{R}}:=N^1(X)\otimes \mathbb{R}.$ Let $\phi:X\rightarrow Y$ be a morphism of smooth projective surfaces. It induces a natural map
$\phi^*: N^1(Y)_{\mathbb{R}}\rightarrow N^1(X)_{\mathbb{R}}$. Since $\dim X=2$, one has a perfect pairing $$N^1(X)_{\mathbb{R}}\times
N^1(X)_{\mathbb{R}}\rightarrow \mathbb{R},\,\,\,
(\delta,\gamma)\rightarrow (\delta\cdot\gamma)\in \mathbb{R}$$ induced by the intersection form. We denote by $\phi_*:N^1(X)_{\mathbb{R}}\rightarrow N^1(Y)_{\mathbb{R}}$ the dual operator of $\phi^*$.

\smallskip

Let $X,Y$ be two smooth projective surfaces and $f:X\dashrightarrow Y$ be a birational map. We denote by $I(f)\subseteq X$ the indeterminacy set of $f$. For any curve $C\subset X$, we write
$$f(C) :=  \overline{ f (C \setminus I(f) ) }$$ the strict transform of $C$.

\smallskip

Let $f:X\dashrightarrow X$ be a birational transformation and $\Gamma$ be a desingularization of its graph. Denote by $\pi_1: \Gamma\rightarrow X,\pi_2: \Gamma\rightarrow X$ the natural projections.
Then the diagram
$$\xymatrix{
& \Gamma \ar[dl]_{\pi_1}\ar[dr]^{\pi_2} &  \\
                     X \ar@{-->}[rr]^f    &  & X
}  \eqno (*) $$

is commutative and we call it a $resolution$ of $f$.

\begin{pro}[\cite{Hartshorne1977}]\label{resolution}
We have the following properties.
\begin{points}
\item
The morphisms $\pi_1,\pi_2$ are compositions of point blowups.
\item
For any point $Q\not\in I(f)$, there is a Zariski open neighborhood $U$ of $p$ in $X$ and an injective morphism $\sigma: U\rightarrow \Gamma$ such that $\pi_1\circ\sigma=\id.$
\end{points}
\end{pro}

Then we define the following linear maps $$f^*=\pi_{1*} \pi_2^*: N^{1}(X)_{\mathbb{R}}\rightarrow N^{1}(X)_{\mathbb{R}},$$ and $$f_*=\pi_{2*} \pi_1^*: N^{1}(X)_{\mathbb{R}}\rightarrow N^{1}(X)_{\mathbb{R}}.$$ Observe that $f_*=f^{-1*}$.
Note that in general we have $(f\circ g)^*\neq g^*f^*.$

\medskip

For any big and nef class $\omega \in N^1_\R(X)$, we set
$$ \deg_\omega(f):=(f^*\omega\cdot \omega),$$
the limit
$\lim_{n \rightarrow\infty}\deg_{\omega}(f^n)^{1/n}$ exists and does not depend on the choice of $\omega$ (see \cite{favre,Dinh2005}). We denote this limit by $\lambda(f)$ and call it the $dynamical$ $degree$ of $f.$

\begin{defi}[see \cite{favre}]\label{defas}
Let $f:X\dashrightarrow X$ be a birational transformation on a smooth projective surface. Then $f$ is said to be $algebraically$ $stable$ if and only if
there is no curve $V\subseteq X$ such that $f^n(V)\subseteq I(f)$ for some integer $n\geq 0.$
\end{defi}
In the case $X=\mathbb{P}^2$, $f$ is algebraically stable if and only if $\deg (f^n)=(\deg f)^n$ for any $n\in \mathbb{N}.$

\begin{thm}[\cite{favre}]\label{AS modele}
Let $f:X\dashrightarrow X$ be a birational transformation of a
smooth projective surface. Then there exists a smooth projective surface $\widehat{X}$, and a proper modification $\pi:
\widehat{X}\rightarrow X$ such that the lift of $f$ to $\widehat{X}$ is an algebraically stable map.
\end{thm}
\bigskip

By a $compactification$ of $\mathbb{A}^2$, we mean
a smooth projective surface
$X$ admitting a birational morphism $\pi : X\dashrightarrow \mathbb{P}^2$
that is an isomorphism above $\mathbb{A}^2\subseteq \mathbb{P}^2$, see \cite{Favre2011}.

The theorem follows from \cite[Proposition 2.6]{Favre2011} and \cite[Theorem 3.1]{Favre2011}, and provides us with a good compactification of $\mathbb{A}^2$.

\begin{thm}[\cite{Favre2011}]\label{l>1goodmodel}Let $f:\mathbb{A}^2\rightarrow \mathbb{A}^2$
be a birational polynomial transformation with $\la (f)>1$. Then there exists a compactification $X$ of $\mathbb{A}^2$ satisfying the following properties.
\begin{points}
\item
The map $f$ extends to an algebraically stable map $\widetilde{f}$ on $X$.
\item
There exists a $\widetilde{f}$-fixed point $Q\in X \setminus \mathbb{A}^2$ such that $d\widetilde{f}^2(Q)=0$.
\item
There exists an integer $n\geq 1$ such that $\widetilde{f}^n(X \setminus \mathbb{A}^2)=Q.$
\end{points}
\end{thm}

\subsection{Branches of curves on surfaces}\cite{Fulton1969,Kunz2005}\label{subsecgerm}
Let $X$ be a smooth projective surface over an algebraically closed field $k$. Let $C$ be an irreducible curve in $X$ and $p$ be a point in $C.$

\begin{defi}A $branch$ of $C$ at $p$ is a point in the normalization of $C$ whose image is $p.$
\end{defi}
\medskip

Let $I_{C,p}$ be the prime ideal associated to $C$ in the local function ring $O_{X,p}$ at $p$ and $\widehat{I_{C,p}}$ be the completion of $I_{C,p}$ in the completion of local function ring $\widehat{O_{X,p}}$.

Let $i:\widetilde{C}\rightarrow C$ is a normalization of $C$ and $\widetilde{p}$ a point in $i^{-1}(p)$. Let $s$ be the branch of $C$ at $p$ defined by the point $\widetilde{p}.$
The morphism $i:\widetilde{C}\rightarrow C$ induces a morphism $i^*:\widehat{\mathcal{O}_{X,p}}\rightarrow\widehat{\mathcal{O}_{\widetilde{C},\widetilde{p}}}$ between the completions of local function rings.
The map $s\mapsto \mathfrak{p}_{s}:=\ker i^*$ gives us a one to one correspondence between the set of branches of $C$ at $p$ and the set of prime ideals of $\widehat{O_{X,p}}$ with height $1$ which contains $\widehat{I_{C,p}}$.

Given any two different branches $s_1$ and $s_2$ at a point $p\in X$, the intersection number is denoted by
$$(s_1\cdot s_2):=
\dim_{k}\widehat{\mathcal{O}_{X,p}}/(\mathfrak{p}_{s_1}+\mathfrak{p}_{s_2}).
$$
For convenience, we set $(s_1\cdot s_2):=0$ if $s_1$ and $s_2$ are branches at different points.

\medskip

Let $Z$ be a smooth projective surface and $f:X\dashrightarrow Z$ be a birational map. If $f$ does not contract $C$ then we denote by $f(s)$ the branch of $f(s)$ defined by the point $\widetilde{p}$ in the normalization $f\circ i: \widetilde{C}\rightarrow f(C)$ and call it the strict transform of $s$. Observe that $f(s)$ is a branch of $f(C)$ and when $p\not\in I(f)$, we have that $f(s)$ is a branch of curve at $f(p).$

If $f$ is regular at $p$, we write
$$f_*s=\left\{ \begin{array}{ccccccccc}
f(s),\text{ when }f \text{ does not contract } C; \\
0,\text{ otherwise}.
\end{array}\right.$$

Let $Y$ be another smooth projective surface and $\pi: Y\rightarrow X$ be a birational morphism.  Denote by $\pi^{\#} s:=\pi^{-1}(s)$ the strict transform of $s.$ Let $E_i$, $i=1,\cdots ,m$ be the exceptional curves of $\pi$. There is a unique sequence of non negative integers $(a_i)_{0\leq i\leq m}$ such that for any irreducible curve $D$ in $Y$ different from $\pi^{\#} C$, we have $(s\cdot \pi_*D)=(\pi^{\#} s+\sum_{i=1}^ma_iE_i\cdot D)$. Denote by $\pi^*s:=\pi^{\#} s+\sum_{i=1}^ma_iE_i$ and call it the pull back of $s.$

\begin{pro}We have the following properties.
\begin{points}
\item We have $\pi_*\pi^*s=s.$
\item For any irreducible curve (resp. any branch of curve) $D$ in $Y$ different from $\pi^{\#}C$ (resp. $\pi^{\#}s$), we have $$(\pi^*s\cdot D)=(s\cdot \pi_*D).$$
\item For any curve (resp. any branch of curve) $D$ in $X$ different from $C$ (resp. $s$) then we have $$(s\cdot D)=(\pi^{\#}s\cdot \pi^*D).$$

\end{points}
\end{pro}

\subsection{Metrics on projective varieties defined over a valued field}\label{secmop}

A field with an absolute value is called a valued  field.
\medskip

\begin{defi}Let $(K,|\cdot|_v)$ be a valued field. For any integer $n\geq 1,$ we define a metric $d_v$ on the projective space $\mathbb{P}^n(K)$ by
$$d_v([x_0:\cdots:x_n],[y_0:\cdots:y_n])=\frac{\max_{0\leq i,j\leq n}|x_iy_j-x_jy_i|_v}{\max_{0\leq i\leq n}|x_i|_v\max_{0\leq j\leq n}|y_j|_v}$$ for any two points $[x_0:\cdots:x_n],[y_0:\cdots:y_n]\in \mathbb{P}^n(K).$
\end{defi}
\smallskip

Observe that when $|\cdot|_v$ is archimedean,  then the metric $d_v$ is not induced by a smooth riemannian metric. However it is equivalent to the restriction of the Fubini-Study metric on $\mathbb{P}^n(\mathbb{C})$ or $\mathbb{P}^n(\mathbb{R})$ to $\mathbb{P}^n(K)$ induced by $\sigma_v$.

\bigskip

More generally, for a projective variety $X$ defined over $K$, if we fix an embedding $\iota: X\hookrightarrow\mathbb{P}^n$, we may restrict the metric $d_v$ on $\mathbb{P}^n(K)$ to a metric $d_{v,\iota}$ on $X(K).$ This metric depends on the choice of embedding $\iota$  in general, but for different embeddings $\iota_1$ and $\iota_2$, the metrics $d_{v,\iota_1}$ and $d_{v,\iota_2}$ are equivalent.
Since we are mostly intersecting in the topology induced by these metrics we shall usually write $d_v$ instead of $d_{v,\iota}$ for simplicity.

\section{A criterion for a curve to be periodic}\label{acf}
Our aim in this section is to prove Theorem \ref{qinf^nC} from the introduction.
Let us recall the setting: \begin{points}
\item
 $X$ is a smooth projective surface over an algebraically closed field;
\item
 $f:X\dashrightarrow X$ is an algebraically stable birational transformation;
\item
$C$ is an irreducible curve in $X$ such that $f^n$ does not contract $C$ and $f^n(C)\bigcap I(f)\neq \emptyset$ for all $n$.
\end{points}
Our aim is to show that $C$ periodic.
\bigskip
Let us begin with the following special case.
\begin{lem}\label{simplercase}Let $x$ be a point in $I(f)\bigcap C$. If there exists a branch $s$ of $C$ at $x$ such that $f^n(s)$ is again a branch at $x$ for all $n\geq 0$, then $C$ is fixed by $f$.
\end{lem}
\proof[Proof of Lemma \ref{simplercase}]
Since $f$ is birational, we may chose a resolution of $f$ as in the diagram $(*)$ in Section \ref{bobmos}.

If $C$ is not fixed, we have $f(s)\neq s$ so that $A:=(s\cdot f(s))_x<\infty$.  By Proposition \ref{resolution}, $\pi_2$ is invertible on a Zariski neighbourhood of $x$. Let $F_x$ be the fiber of $\pi_1$ over $x.$

 For any $m\geq 0,$ we have,
 \begin{align*}
((f^{m}(s)\cdot f^{m+1}(s))_x=&\sum_{y\in F_x}(\pi_1^{\#}f^{m}(s)\cdot \pi_1^*f^{m+1}(s))_y\nonumber\\
        \geq&(\pi_1^{\#}f^{m}(s)\cdot \pi_1^{*}f^{m+1}(s))_{\pi_2^{-1}(x)}\\
         =&(\pi_1^{\#}f^{m}(s)\cdot \pi_1^{\#}f^{m+1}(s))_{\pi_2^{-1}(x)}+(\pi_1^{\#}f^{m}(s)\cdot F_x)_{\pi_2^{-1}(x)}\\
         =&(f^{m+1}(s)\cdot f^{m+2}(s))_{x}+(\pi_1^{\#}f^{m}(s)\cdot F_x)_{\pi_2^{-1}(x)}\\
         \geq& (f^{m+1}(s)\cdot f^{m+2}(s))_{x}+1.
\end{align*}
 It follows that $A=(s\cdot f(s))_x\geq (f^{m}(s)\cdot f^{m+1}(s))_x+m\geq m$ for all $m\geq 0$ which yields a contradiction.
\endproof

We now treat the general case.
\proof[Proof of Theorem \ref{qinf^nC}]Recall that $f^n$ does not contract $C$ and $f^n(C)\bigcap I(f)\neq \emptyset$ for all $n$. By Lemma \ref{simplercase}, it is sufficient to find a point $x\in I(f)\bigcap C$ such that the image by $f^n$ of the branch of $C$ at $x$ is again a branch of a curve at $x$ for all $n\geq 0.$ By contradiction we suppose that $C$ is not periodic.

To do so, we introduce the set $$P(f)=\{x\in I(f)|\,\, \text{ there is }n_1>n_2\geq 0 \text{ such that } f^{-n_1}(x)=f^{-n_2}(x)\}$$ and the set $$O(f)=\{f^{-n}(x)|\,\,x\in P(f) \text{ and }n\geq 0\}.$$ By definition, $O(f)$ is finite. Since $f$ is algebraically stable, $O(f)=O(f^n)$ for all $n\geq 1.$ Replacing $f$ by $f^l$ for a suitable $l\geq 1$,  we may assume that $O(f)=P(f).$  Set $N(f)=I(f)\setminus P(f).$

\medskip
First, we prove
\begin{lem}\label{iftoof} For all $n\geq 0$, $f^n(C)\bigcap O(f)\neq\emptyset.$
\end{lem}
\proof[Proof of Lemma \ref{iftoof}]
 We assume that $I(f)=\{p_1,\cdots, p_m\}$ and define the map $$F=(f^{-1},\cdots,f^{-1}):X^m\dashrightarrow X^m.$$  Denote by $\pi_i$ the projection onto the $i$-th factor and set $$D=\bigcup_{i=1}^m \pi_i^{-1}(C).$$ Pick a point $q=(p_1,\cdots,p_m)\in X^m.$ Since $f^n(C)\bigcap I(f)\neq \emptyset$ for all $n\geq 0$ by assamption, we have $F^n(q)\in D$ for all $n\geq 0.$  Let $Z'$ be the Zariski closure of $\{F^n(q)|n\geq 0\}.$ Then we have $Z'\subseteq D.$ Let $Z$ be the union of all irreducible components of $Z'$ of positively dimension. If $Z$ is empty, then $p_i$ is $f^{-1}$-preperiodic for all $i$ and we conclude.

 Otherwise since $\{F^n(q)|n\geq 0\}\bigcap I(F)=\emptyset,$ the proper transformation of $Z$ by $F$ is well defined and satisfies $F(Z)=Z$, hence all irreducible components of $Z$ are periodic. Let $l$ be a common period for all components of $Z$. Observe that any irreducible component of $Z$ is included in some $\pi^{-1}_i(C)$ for $i=1,\cdots,m.$
 In other words, there exists $k\geq 0$ and $i\in \{1,\cdots,m\}$ such that $f^{-ln-k}(p_i)\in C$ for all $n\geq 0$. If $p_i$ is not $f^{-1}$-preperiodic, then $C$ is the Zariski closure of $\{f^{-ln-k}(p_i)|n\geq 0\}$ which is $f^{-l}$-invariant. This implies $C$ to be periodic which contradicts to our hypothesis.  It follows that $p_i$ is $f^{-1}$-preperiodic.

Repeating the same argument for $f^n(C)$, we have $f^n(C)\bigcap O(f)\neq \emptyset$ for all $n\geq 0.$
\endproof
\medskip
Denote by $D(n)$ the number of branches of $f^n(C)$ at points of $O(f)$.
 Since $f^{-1}(O(f))\subseteq O(f)$, we have $D(n)$ is decrease and by Lemma \ref{iftoof}, we have $D(n)\geq 1.$
Replace $C$ by $f^M(C)$ for some $M\geq 0$, we may assume that $D(n)$ is constant for $n\geq 0.$
It follows that for any branch of curve of $f^n(C)$ at a point in $O(f)$, its image by $f$ is again a branch of $f^{n+1}(C)$ at a point of $O(f)$.
Set $$S=\{x\in O(f)| \text{ there are infinitely many }n\geq 0 \text{ such that }x\in f^n(C)\}.$$

By the finiteness of $O(f)$, we may suppose that $$f^n(C)\bigcap O(f)=f^n(C)\bigcap S$$ for all integer $n\geq 0.$
\medskip

We claim that
 \begin{lem}\label{germfix}Replacing $f$ by a positive iterate, there exists a point $x\in C\bigcap S$ for which there is a branch $s$ of $C$ at $x$ such that $f^n(s)$ is again a branch of curve at $x$ for all $n\geq 0.$
\end{lem}
According to Lemma \ref{simplercase}, we conclude.
\endproof

\proof[Proof of Lemma \ref{germfix}]
Pick a resolution of $f$ as in the diagram $(*)$ in Section \ref{bobmos}.
For any point $x\in S$, denote by $F_x$ the fibre of $\pi_1$ over $x$ and $E_x=\pi_2(F_x)\bigcap S$.
\smallskip

We have $E_x\neq\emptyset.$ Otherwise,
there exists $n\geq 0$ for which $x\in f^n(C)$ and a branch $s$ of $f^n(C)$ at $x$.  The assumption $E_x=\emptyset$ implies that $f(s)$ is not a branch at any point in $S$. This shows that $D(n+1)<D(n)$ and we get a contradiction.

On the other hand, let $x_1,x_2$ be two different points in $S$. If $E_{x_1}\bigcap E_{x_1}\neq \emptyset$, there exists $y\in S$ such that $y\in \pi_2(F_{x_1})\bigcap \pi_2(F_{x_2})$. By Zariski's main theorem, $\pi_2^{-1}(y)$ is a connected curve meeting $F_{x_1}$ and $F_{x_2}$. So $\pi_1(\pi_2^{-1}(y))$ is a curve and it is contracted by $f$ to $y\in S\subseteq I(f)$. This contradicts the fact that $f$ is algebraically stable. So we have $$E_{x_1}\bigcap E_{x_1}=\pi_2(F_{x_1})\bigcap \pi_2(F_{x_2})\bigcap S=\emptyset.$$

Set $T=\coprod_{x\in S} E_x\subseteq S.$ Since $\# E_x\geq 1$ for all $x$, we have $\#T\geq \# S.$ It follows that $T=S$ and $\# E_x=1$ for all $x\in S.$ This allows us to define a map $G:S\rightarrow S$ sending $x\in S$ to the unique point in $E_x$. Then $G$ is an one to one map. For all $n\geq M,$ $f$ sends a branch of $f^n(C)$ at a point $x\in S$ to a branch of $f^{n+1}(C)$ at the point $G(x)$. By replacing $f$ by $f^{(\# S)!}$, we may assume that $G=\id.$  Then for any $x\in S\bigcap C$ and $s$ a branch of $C$ at $x$, we have $f^n(s)$ is again a branch at $x$ for all $n\geq 0.$
\endproof

\section{The DML property}\label{sedml}
For convenience, we introduce the following
\begin{defi}Let $X$ be a smooth surface defined over an algebraically closed field, and $f:X\dashrightarrow X$ be a rational transformation. We say that the pair $(X,f)$ satisfies the DML property if
for any curve $C$ on $X$ and for any closed point $p\in X$ such that $f^n(p)\not\in I(f)$ for all $n\geq 0$, the set $\{n\in \mathbb{N}|\,\, f^n(p)\in C\}$ is a union of at most finitely many arithmetic progressions.
\end{defi}
In our setting the DML property is equivalent to the following seemingly stronger property.
\begin{pro}\label{prodefidml}Let $X$ be a smooth surface defined over an algebraically closed field, and $f:X\dashrightarrow X$ be a rational transformation.
The following statements are equivalent.
\begin{itemize}
\item[(1)] The pair $(X,f)$ satisfies the DML property.
\item[(2)] For any curve $C$ on $X$ and any closed point $p\in X$ such that $f^n(p)\not\in I(f)$ for all $n\geq 0$ and the set $\{n\in \mathbb{N}| f^n(p)\in C\}$ is infinite, then $p$ is preperiodic or $C$ is periodic.
\end{itemize}
\end{pro}

\proof Suppose (1) holds. Let $C$ be any curve in $X$ and $p$ be a closed point in $X$ such that $f^n(p)\not\in I(f)$ for all $n\geq 0$. Assume that the set $\{n\in \mathbb{N}|\,\, f^n(p)\in C\}$ is infinite. The DML property of $(X,f)$ implies that there are integers $a>0$ and $b\geq 0$ such that $f^{an+b}(p)\in C$ for all $n\geq 0.$ If $p$ is not preperiodic, the set $O_{a,b}:=\{f^{an+b}(p)|\,\,n\geq 0\}$ is Zariski dense in $C$ and $f^{a}(O_{a,b})\subseteq O_{a,b}$. It follows that $f^a(C)\subseteq C$, hence $C$ is periodic.

Suppose (2) holds. If the set $S:=\{n\in \mathbb{N}|\,\, f^n(p)\in C\}$ is finite or $p$ is preperiodic, then there is nothing to prove. We may assume that $S$ is infinite and $p$ is not preperiodic. The property (2) implies that $C$ is periodic. There exists an integer $a>0$ such that $f^a(C)\subseteq C$. We may suppose that $f^{i}(C)\not\subseteq C$ for $1\leq i\leq a-1$. Since $p$ is not preperiodic, there exists $N\geq 0$, such that $f^n(p)\not\in (\bigcup_{1\leq i\leq a-1}f^{i}(C))\bigcap C$ for all $n\geq N.$
So $S\setminus\{1,\cdots,N-1\}$ takes form $\{an+b|\,\,n\geq 0\}$ where $b\geq 0$ is an integer, and it follows that $(X,f)$ satisfies the DML property.
\endproof

\begin{thm}\label{ftof^n}\label{open}\label{blowdown}\label{blowup}Let $X$ be a smooth surface defined over an algebraically closed field, and $f:X\dashrightarrow X$ be a rational transformation, then the following properties hold.
\begin{points}
\item For any $m\geq 1$, $(X,f)$ satisfies the DML property if and only if $(X,f^m)$ satisfies the DML property.
\item Suppose $U$ is an open subset of $X$ such that the restriction $f_{|U}:U\rightarrow U$ is a morphism. Then $(X,f)$ satisfies the DML property, if and only if $(U,f_{|U})$ satisfies the DML property.
\item Suppose $\pi:X\rightarrow X'$ is a birational morphism between smooth projective surfaces, and $f:X\dashrightarrow X$, $f':X'\dashrightarrow X'$ are rational maps such that $\pi\circ f=f'\circ\pi$. If the pair $(X,f)$ satisfies the DML property, then $(X',f')$ satisfies the DML property.
\item Suppose $\pi:X\rightarrow X'$ is a birational morphism between smooth projective surfaces, and $f:X\dashrightarrow X$, $f':X'\dashrightarrow X'$ are \emph{birational} transformations such that $\pi\circ f=f'\circ\pi$. If $f'$ is algebraically stable and the pair $(X',f')$ satisfies the DML property, then $(X,f)$ satisfies the DML property.
\end{points}
\end{thm}

\begin{defi}Let $X$ be a smooth projective surface defined over an algebraically closed field and $f:X\dashrightarrow X$ be a birational transformation. We say that $(X',f')$ is a birational model of $(X,f)$ if there is a birational map $\pi: X'\dashrightarrow X$ such that $$f'=\pi^{-1}\circ f\circ\pi.$$
\end{defi}

\begin{cor}\label{astoall}Let $X$ be a smooth projective surface defined over an algebraically closed field and $f:X\dashrightarrow X$ be an algebraically stable birational transformation such that $(X,f)$ satisfies the DML property. Then all birational models $(X',f')$ of $(X,f)$ satisfy the DML property.
\end{cor}

\proof[Proof of Corollary \ref{astoall}]
Pick $Y$ a desingularization of the graph of $f$ and set $\pi_1,\pi_2$ the projections which make the diagram
$$\xymatrix{
& Y \ar[dl]_{\pi_1}\ar[dr]^{\pi_2} &  \\
                     X \ar@{-->}[rr]^{\phi}    &  & X'
}$$
to be commutative. Since $f$ is algebraically stable, its lift to $Y$ satisfies the DML property by Theorem \ref{blowup} (iv). We conclude that $(X',f')$ satisfies the DML property by Theorem \ref{blowdown} (iii).
\endproof

\proof[Proof of Theorem \ref{ftof^n}]
\proof[(i)] The "only if" part is trivial, so that we only have to deal with the "if" part. We assume that $(X,f^m)$ satisfies the DML property. Let $C$ be a curve in $X$ and $p$ be a point in $X$ such that $f^n(p)\not\in I(f)$ for all $n\geq 0$. Suppose that the set $\{n\in \mathbb{N}|\,\, f^n(p)\in C\}$ is infinite. Since $$\{n\in \mathbb{N}|\,\, f^n(p)\in C\}=\bigcup_{i=0}^{m-1}\{n\in \mathbb{N}|\,\, f^{nm}(f^i(p))\in C\},$$ then for some $i$, the set $\{n\in \mathbb{N}| \,\,f^{nm}(f^i(p))\in C\}$ is also infinite. Since $(X,f^m)$ satisfies the DML property, $C$ is periodic or $f^i(p)$ is preperiodic. It follows that $C$ is periodic or $p$ is preperiodic.
\proof[(ii)] If $(X,f)$ satisfies the DML property, since $f_{|U}:U\rightarrow U$ is a morphism, $(U,f_{|U})$ satisfies the DML property.

Conversely suppose that $(U,f_{|U})$ satisfies the DML property. Let $C$ be an irreducible curve in $X$, $p$ be a closed point in $X$ such that $f^n(p)\not\in I(f)$ for all $n\geq 0$ and the set $\{n\in \mathbb{N}| f^n(p)\in C\}$ is infinite. The set $E=X-U$ is a proper closed subvariety of $X$. If $p\in U$, then we have that $C\not\subseteq E$. Since $(U,f_{|U})$ satisfies the DML property, we have either $p$ is preperiodic or $C$ is periodic. Otherwise, we may assume that for all $n\geq 0,$ $f^n(p)\in E,$ then the Zariski closure $D$ of $\{f^n(p)|\,\,n\geq 0\}$, is contained in $E$. We assume that $p$ is not preperiodic, then $C\subseteq D$. Since $D$ is fixed, we have that $C$ is periodic.

\proof[(iii)] It is sufficient to treat the case when $\pi$ is the blowup at a point $q\in X'.$ Let $C'$ be a curve in $X'$, $p'$ be a point in $X'$ such that $(f')^n(p')\not\in I(f')$ for all $n\geq 0$ and the set $\{n\in \mathbb{N}|\,\, (f')^n(p')\in C'\}$ is infinite. We assume that $p'$ is not a periodic point, so that for $n$ large enough, $f^{'n}(p')\neq q$. Replacing $p$ by $f^{'m}(p')$ for some $m$ large enough, we may assume that $f^{'n}(p')\neq q$ for all $n\geq 0.$ Set $p=\pi^{-1}(p')$ and $C=\pi^{-1}(C')$, then we have $f^n(p)\not\in I(f)$ for all $n\geq 0$ and the set $$\{n\in \mathbb{N}|\,\, f^n(p)\in C\}$$ is infinite. This implies $C$ and then $C'$ to be periodic.
\proof[(iv)] Let $C\subseteq X$ be a curve, $p$ be a point in $X$ such that $f{^n}(p)\not\in I(f)$ for all $n\geq 0$ and the set $\{n\in \mathbb{N}|\,\, f{^n}(p)\in C\}$ is infinite. We may assume that $C$ is irreducible. Let $E$ be the exceptional locus of $\pi.$

\begin{lem}\label{ccEqcI}If $C\subseteq E$ and $\pi(C)$ is a point in $I(f')$, then (iv) holds.
\end{lem}
\proof[Proof of Lemma \ref{ccEqcI}]Set $q:=\pi(C)\in I(f')$. Since $f'$ is algebraically stable, we have $q\not\in I((f^{'})^{-n})$ and $$\pi(f^{-n}(C))=(f^{'})^{-n}(q)$$ for all $n\geq 1.$ It follows that $f^{-n}(C)$ is a point or an exceptional curve of $\pi$ for $n\geq 1.$

If there exists $l\geq 1$ such that $f^{-l}(C)$ is a point, we pick two integers $n_1>n_2\geq l$ such that $f^{n_1}(p),f^{n_2}(p)\in C.$
 Then $f^{n_1-l}(p)=f^{n_2-l}(p)$, which implies $p$ to be preperiodic.

  Otherwise $f^{-n}(C)$ is an exceptional curve of $\pi,$ for all $n\geq 0$. Since there are only finitely many irreducible components of $E$, we have that $C$ is periodic.
\endproof

Denote by $K=\pi^{-1}(I(f')).$
\begin{lem}\label{infK}If there are infinitely many $n\geq 0$ such that $f^n(p)\in K,$ then (iv) holds.
\end{lem}
\proof[Proof of \ref{infK}]
There is an irreducible component $F$ of $K$ such that the set $\{n\geq0 |\,\,f^n(p)\in F\}$ is infinite.

If $F$ is a point, then $p$ is preperiodic.

Otherwise $F$ is a curve, then $F\subseteq E$ and $\pi(F)\subseteq I(f').$ Suppose that $p$ is not preperiodic, Lemma \ref{ccEqcI} shows that $F$ is periodic. Then $F'=\bigcup_{k\geq 0}f^k(F)$ is a curve and $f^{n}(p)\subseteq F'$ for all $n\geq 0.$ If $C\subseteq F'$, then $C$ is periodic. If $C\not\subseteq F',$ then $C\bigcap F'$ is finite, and this shows that $p$ is preperiodic.
\endproof

\begin{lem}\label{ccE} If $C\subseteq E$, then (iv) holds.
\end{lem}

\proof
By Lemma \ref{infK}, we may assume that there exists an integer $N\geq 0$, such that $f^n(p)\not\in K$ for all $n\geq N.$

Set $q:= \pi(C)$. By Lemma \ref{ccEqcI}, we assume that $q\not\in I(f')$. Then we have $$\pi(f^{N+l}(p))=f^{'l}(\pi(f^N(p)))$$ for $l\geq 0$. It follows that there are infinitely many $l\geq 0$, such that $f^{'l}(\pi(f^N(p)))=q.$ Then $q$ is preperiodic and the obit of $f^{'N}(q)$ does not meet $I(f')$. Since $\pi(f^{n}(C))=f^{'n}(q)$ for all $n\geq 0$, we have $f^n(C)\subseteq \bigcup_{k\geq N}\pi^{-1}(f^k(q))$ for all $n\geq N.$ Hence ether $C$ is periodic or for some $n\geq 1,$ $f^n(C)$ is a point. In the second case, we conclude that $p$ is preperiodic.
\endproof
Let $L=K\bigcup E.$
\begin{lem}\label{infL}If there are infinitely many $n\geq 0$ such that $f^n(p)\in L,$ then (iv) holds.
\end{lem}
\proof[Proof of Lemma \ref{infL}]
There is an irreducible component $F$ of $L$ such that $\{n\geq0 |\,\,f^n(p)\in F\}$ is infinite.

If $F$ is a point, then $p$ is preperiodic.

Otherwise $F$ is a curve, then $F\subseteq E.$ Suppose that $p$ is not preperiodic, Lemma \ref{ccE} shows that $F$ is periodic. Then $F'=\bigcup_{k\geq 0}f^k(F)$ is a curve and $f^{n}(p)\subseteq F'$ for all $n\geq 0.$ If $C\subseteq F'$, then $C$ is periodic. Otherwise $C\not\subseteq F',$  we have that $C\bigcap F'$ is finite and then $p$ is preperiodic.
\endproof
We may assume that there is an integer $M\geq 0$, such that $f^n(p)\not\in L$ for all $n\geq M.$

If $C\not\subseteq E$, $\pi(C)$ is a curve. For all $l\geq 0$ we have $$\pi(f^{M+l}(p))=f^{'l}(\pi(f^M(p)))\not\in I(f').$$ Since $(X',f')$ satisfies the DML property, either $\pi(C)$ is periodic or $\pi(p)$ is preperiodic. When $\pi(C)$ is periodic, we have $C$ is periodic. Otherwise $\pi(p)$ is preperiodic. For any $l\geq 0$, $\pi$ is invertible on some Zariski neighborhood of the point $f^{'l}(\pi(f^M(p)))$ and then we conclude that $p$ is peperiodic.
\endproof

\section{Local dynamical Mordell Lang theorem}\label{seldml}
The aim of this section is to prove Theorem \ref{local}.
We are in the following situation:
\begin{points}
\item $X$ is a smooth projective surface defined over an arbitrary valued field $(K,|\cdot|).$
\item $f:X\dashrightarrow X$ is a birational transformation defined over $K;$
\item $Q$ is $K$-point of $X$ such that $Q\in I(f^{-1})\bigcap I(f)$ and $f(Q)=Q;$
\item $p$ is $K$-point of $X$ such that $f^n(p)\not\in I(f)$ for all $n\geq 0;$
\item $f^n(p)\rightarrow Q$ as $n\rightarrow\infty$ with respect to the topology induced by $|\cdot|;$
\item $C$ is a curve in $X$ such that the set $\{n\in \mathbb{N}|\,\,f^n(p)\in C\}$ is infinite.
\end{points}
We want to prove that $C$ is fixed by $f.$

\proof[Proof of Theorem \ref{local}]
Pick a resolution of $f$ as in the diagram $(*)$ in Section \ref{bobmos}. Recall Proposition \ref{resolution}.
Assume that for all $n\geq 0$, $f^n(p)\neq Q.$ There is an infinite sequence $\{n_k\}_{k\geq 0}$ such that $f^{n_k}(p)\in C\setminus\{Q\}$. It follows that $f^{n_k-m}(p)\in f^{-m}(C)$ for $k$ large enough. Setting $k\rightarrow \infty$, we get $Q\in f^{-m}(C)$ for all $m\geq 0.$

If $C\neq f^{-1}(C)$, then we have $f^{-m}(C)\neq f^{-m-1}(C)$ for all $m\geq 0$. By computing local intersection at $Q$, we get
\begin{multline}
(f^{-m}(C)\cdot f^{-m-1}(C))_Q=\sum_{x\in \pi_2^{-1}(Q)}(\pi_2^*f^{-m}(C)\cdot \pi_2^{\#}f^{-m-1}(C))_x\\
         =\sum_{x\in \pi_2^{-1}(Q)}\left(\left(\pi_2^{\#}f^{-m}(C)+\sum_{i=1}^sv_{E_i}(f^{-m}(C))E_i\right)\cdot \pi_2^{\#}f^{-m-1}(C)\right)_x
         \end{multline}
where $E_i, 1\leq i \leq s$ are irreducible exceptional curves for $\pi_2.$ Since $$\Supp (\sum_{i=1}^sv_{E_i}(f^{-m}(C))E_i)=\bigcup_{1\leq i \leq s}E_i=\pi_2^{-1}(Q),$$ we have

\begin{multline*}
(5.1)=\sum_{x\in \pi_2^{-1}(Q)}(\pi_2^{\#}f^{-m}(C)\cdot \pi_2^{\#}f^{-m-1}(C))_x+\left(\left(\sum_{i=1}^sv_{E_i}(f^{-m}(C))E_i\right)\cdot \pi_2^{\#}f^{-m-1}(C)\right)\\
\geq\sum_{x\in \pi_2^{-1}(Q)}(\pi_2^{\#}f^{-m}(C)\cdot \pi_2^{\#}f^{-m-1}(C))_x+1\\
=\left(\sigma (f^{-m-1}(C))\cdot \sigma(f^{-m-2}(C))\right)_{\sigma (Q)}+1\\
=(f^{-m-1}(C)\cdot f^{-m-2}(C))_{Q}+1.
        \end{multline*}

It follows that $$0<(f^{-m}(C)\cdot f^{-m-1}(C))_{Q}\leq (f^{-m+1}(C)\cdot f^{-m}(C))_{Q}-1\leq \cdots \leq (C\cdot f^{-1}(C))_{Q}-m$$ for all $m\geq 0,$ which yields a contradiction. So we have $C=f^{-1}(C)$ and then $f(C)=C.$
\endproof

Observe that our proof of Theorem \ref{local} actually gives
\begin{pro}\label{periodictofixed}Let $X$ be a projective surface over an algebraically closed field and $f: X\dashrightarrow X$ be a birational map with a fixed point $Q\in I(f^{-1})\setminus I(f)$. Then all periodic curves passing through $Q$ are fixed.
\end{pro}
\section{The case $\lambda(f)=1$}\label{l=1}
In this section, we prove Theorem A in the case $\lambda(f)=1$. Denote by $K$ an algebraically closed field of characteristic $0$.

Recall from \cite{favre} and \cite{Favre2011}, that if $\lambda(f)=1$, then we are in one of the following
two cases:
\begin{enumerate}
\item there exists a smooth projective surface $X$ and an automorphism $f'$ on $X$ such that the pair $(X,f')$ is birationally conjugated to $(\mathbb{A}^2,f)$;
\item in suitable affine coordinates, $f(x,y)=(ax+b,A(x)y+B(x))$ where $A$ and $B$ are polynomials with $A\neq 0$ and $a\in K^*$, $b\in K$.
\end{enumerate}

The case of automorphism has been treated by Bell, Ghioca and Tucker. Theorem A thus follows from \cite[Theorem 1.3]{Bell2010} in case (1) and in case (2) where $\deg A=0$.
So
in this section we suppose that $f$ takes form $$f(x,y)=(ax+b,A(x)y+B(x))  \eqno(**)$$ with $A,B\in K[x]$, $\deg A\geq 1 $, $a\in K^*$ and $b\in K$.
\subsection{Algebraically stable models}\label{l=1model}
Any map of the form $(**)$ can be made algebraically stable in a suitable Hirzebruch surface $\mathbb{F}_n$ for some $n\geq 0.$ It is convenient to work with the presentation of these surfaces as a quotient by $(\mathbb{G}_m)^2,$ as in \cite{AntonioLaface}. By definition, the set of closed points $\mathbb{F}_n(K)$ is the quotient of $\mathbb{A}^4(K)\setminus (\{x_1=0 \text{ and } x_2=0\}\bigcup \{x_3=0 \text{ and } x_4=0\})$ by the equivalence relation generated by $$(x_1,x_2,x_3,x_4)\sim (\lambda x_1,\lambda x_2,\mu x_3,\mu/\lambda^n x_4)$$ for $\lambda,\mu \in K^*$. We denote by $[x_1,x_2,x_3,x_4]$ the equivalence class of $(x_1,x_2,x_3,x_4)$. We have a natural morphism
$\pi_n:\mathbb{F}_n\rightarrow \mathbb{P}^1$ given by $\pi_n([x_1,x_2,x_3,x_4])=[x_1:x_2]$ which makes $\mathbb{F}_n$ into a locally trivial $\mathbb{P}^1$ fibration.

We shall look at the embedding $$i_n:\mathbb{A}^2\hookrightarrow \mathbb{F}_n:(x,y)\mapsto [x,1,y,1].$$ Then $\mathbb{F}_n\setminus \mathbb{A}^2$ is union of two lines: one is the fiber at infinity $F_{\infty}$ of $\pi_n$, and the other one is a section of $\pi_n$ which we denote by $L_{\infty}$.
\bigskip

Recall that $f$ has the form $(**)$. For each $n\geq 0,$ set $d=\max\{\deg A, \deg B-n\}$. By the embedding $i_n$, the map $f$ extends to a birational transformation $$f_n:[x_1,x_2,x_3,x_4]\mapsto [ax_1+bx_2,x_2,A(x_1/x_2)x_2^dx_3+B(x_1/x_2)x_2^{d+n}x_4,x_2^{d}x_4]$$ on $\mathbb{F}_n$.
For any $n\geq \deg B-\deg A+1$, we have $d=\deg A$ and $$I(f_n)=\{[x_1,x_2,x_3,x_4]\in \mathbb{F}_n| x_2=x_3=0\}.$$ The unique curve which is contracted by $f_n$ is $F_{\infty}=\{x_2=0\}$ and its image is $f_n(F_{\infty})=[1,0,1,0]$. It implies the following:

\begin{pro}\label{propolycontra}\label{corpolyem}For any integer $n\geq \deg B-\deg A+1$, $f_n$ is algebraically stable on $\mathbb{F}_n$ and contracts the curve $F_{\infty}$ to the point $[1,0,1,0]$.
\end{pro}

\subsection{The attracting case}In the remaining of this section, we fix an integer $m$ such that the extension of $f$ to $\mathbb{F}_m$ is algebraically stable. For simplicity, we write $f$ for the map $f_m$ induced by $f$ on $\mathbb{F}_m$.
\medskip

\begin{pro}\label{A^2}Let $|\cdot|$ be an absolute value on $K$ such that $|a|>1$. Then $(\mathbb{F}_m,f)$ satisfies the DML property.
\end{pro}
\proof Since $a\neq 1$, by changing coordinates, we may assume that $f=(ax,A(x)y+B(x))$. Since $f$ contracts the fiber $F_{\infty}$ to $O:=L_{\infty}\bigcap F_{\infty}$, the point $O$ is fixed and the two eigenvalues of $df$ at $O$ are $1/a$ and $0$.
Since $|a|>1$, there is a neighbourhood $U$ of $O$, such that $U\bigcap I(f)=\emptyset$, $f(U)\subseteq U$ and $f^n\rightarrow O$ uniformly on $U.$

Let $C$ be an irreducible curve in $\mathbb{P}^2_K$ and $p$ be a point in $\mathbb{A}^2_K$ such that the set $\{n\in \mathbb{N}|\,\, f^n(p)\in C\}$ is infinite. By Lemma \ref{open}, we may assume that $p\in \mathbb{A}^2_{K}$ and $C\not\subseteq L_{\infty}\bigcup F_{\infty}.$

If $C\bigcap F_{\infty}=\{O\}$, there is an open set $V$ of $\mathbb{P}^1_{K}$, such that $[1:0]\in V$ and $\pi_m^{-1}(V)\bigcap C\subseteq U.$ Since $|a|>1$, for $n$ large enough, $f^n(p)\in \pi_m^{-1}(V)$. So there is an integer $n_1>0$ such that $f^{n_1}(p)\in U.$ Theorem \ref{local} implies that the curve $C$ is fixed.

\smallskip

We may assume now that $f^n(C)\bigcap F_{\infty}\neq \{O\}$ for all $n\geq 0.$

If $C\bigcap F_{\infty}=\emptyset,$ then $C$ is a fiber of the rational fibration $\pi_m:\mathbb{F}_m\rightarrow \mathbb{P}^1$. Since $\{n\in \mathbb{N}| \,\, f^n(p)\in C\}$ is infinite, the curve $C$ is fixed.

Finally assume that $f^n(C)\bigcap F_{\infty}\neq\emptyset$ for all $n\geq 0.$ Since $f$ contracts $F_{\infty}$ to $O$, we have $$f^n(C)\bigcap I(f)\neq \emptyset,$$ and we conclude by Theorem \ref{qinf^nC} that $C$ is periodic in this case.
\endproof

\subsection{The general case}\label{l=1g}

\begin{pro}\label{l=1as}The pair $(\mathbb{F}_m, f)$ satisfies the DML property.
\end{pro}

\proof Let $C$ be a curve in $\mathbb{F}_m$, and $p$ be a point in $\mathbb{A}_{K}^2$ such that the set $\{n\geq 0| f^n(p)\in C\}$ is infinite.
We may assume that the transcendence degree of $K$ is finite, since we can find a subfield of $K$ such that it has finite transcendence degree and $f, C$ and $p$ are all defined over this subfield.

In the case $f$ acts on the base as the identity, the proposition holds trivially. Assume that it is not that case. Let $O=L_{\infty}\bigcap F_{\infty}$. As in the proof of Proposition \ref{A^2}, we only have to consider the case $C\bigcap F_{\infty}=O.$

\smallskip

If $a$ is a root of unity, we may replace $f$ by $f^n$ for some integer $n>0$ and assume that $a=1$ and $b=1.$ Since the transcendence degree of $K$ is finite, we may embed $K$ in the field of complex numbers $\mathbb{C}.$ Let $|\cdot|$ be the standard absolute value on $\mathbb{C}$.  Since $f$ contracts $F_{\infty}$ to $O$, there is a neighborhood $U$ of $O$ with respect to the usual euclidian topology such that for all point $q\in U\bigcap \{(x,y)\in \mathbb{C}^2|\,\, \text{Re}(x)> 0\}$, we have $\lim_{n\rightarrow \infty} f^n(q)=O.$ Since $C\bigcap F_{\infty}=O,$ there exists $M>0,$ such that $C\bigcap \{(x,y)|\,\, \text{Re}(x)> M\}\subseteq U$ and we conclude by using Theorem \ref{local} in this case.

\smallskip

If $a$ is an algebraic number over $\mathbb{Q}$ and is not a root of unity, by \cite[Theorem 3.8]{Silverman2007} there exists an absolute value $|\cdot|_v$ (either archimedean or
non-archimedean) on $\overline{\mathbb{Q}}$ such that $|a|_v>1.$
This shows that $(\mathbb{F}_m, f)$ satisfies the DML property by Proposition \ref{A^2}.

\smallskip

If $a$ is not an algebraic number over $\mathbb{Q}$, we claim that there exists a field embedding $\iota: K\hookrightarrow \mathbb{C}$ such that $|\iota(a)|>1,$  and we may conclude again by using Proposition \ref{A^2}.

It thus remains to prove the claim. There is a subring $R$ of $K$ which is finitely generated over $\overline{\mathbb{Q}}$, such that $f, C$ and $p$ are all defined over $R.$ There is an integer $l>0$, such that $R=\overline{\mathbb{Q}}[t_1,\cdots, t_l]/I$, where $I$ is a prime ideal of $\overline{\mathbb{Q}}[t_1,\cdots, t_l]$. It induces an embedding $\Spec R:=V\subseteq \mathbb{A}^l_{\overline{\mathbb{Q}}}$. We set $$V_{\mathbb{C}}:=V\times_{\Spec \overline{\mathbb{Q}}}\Spec \mathbb{C} \subseteq \mathbb{A}^l_{\mathbb{C}}.$$ For any polynomial $F\in \overline{\mathbb{Q}}[t_1,\cdots, t_l]\setminus I$, we also define $V_F:=\{F=0\}$. Then $V_{\mathbb{C}}\setminus V_F$ is a dense open set in the usual euclidian topology. Since $\overline{\mathbb{Q}}[t_1,\cdots, t_l]\setminus I$ is countable, the set $V_{\mathbb{C}}\setminus (\bigcup_{F\in \overline{\mathbb{Q}}[t_1,\cdots, t_l]\setminus I}V_F)$ is dense. Interpreting $a$ a nonconstant holomorphic function on $V_{\mathbb{C}}$, we see that there exists an open set $W\subseteq V_{\mathbb{C}}$ such that $|a|>1$ on $W.$

Pick a closed point $(s_1,\cdots, s_l)\in W\setminus (\bigcup_{F\in \overline{\mathbb{Q}}[t_1,\cdots, t_l]\setminus I}V_F)$ and consider the unique morphism $\iota: R=\mathbb{Q}[t_1,\cdots, t_l]/I\rightarrow \mathbb{C}$ sending $t_i$ to $s_i.$ This morphism is in fact an embedding. We may extend it to an embedding of $K$ as required.\endproof
\section{Upper bound on heights when $\lambda(f)>1$}\label{seupheight}
\subsection{Absolute values on fields}(\cite{Silverman2007})
Set $\mathcal{M}_{\mathbb{Q}}:=\{|\cdot|_{\infty} \text { and }|\cdot|_p \text{ for all prime p}\}$
where $|\cdot|_{\infty}$ is the usual absolute value and $|\cdot|_p$ is the $p$-adic absolute value defined by $|x|:=p^{-\ord_p(x)}$ for $x\in \mathbb{Q}$.

\medskip

Let $K/\mathbb{Q}$ be a number field. The $set$ $of$ $places$ on $K$ is denoted by $\mathcal{M}_{K}$ and consists of all absolute values on $K$ whose restriction to $\mathbb{Q}$ is one of the places in $\mathcal{M}_{\mathbb{Q}}.$ Further we denote by $\mathcal{M}_{K}^{\infty}$ the set of archimedean places; and by  $\mathcal{M}_{K}^{0}$ the set of nonarchimedean places.

When $v$ is archimedean, there exists an embedding $\sigma_v:K\hookrightarrow \mathbb{C} \text{ (or } \mathbb{R})$ such that $|\cdot|_v$ is the restriction to $K$ of the usual absolute value on $\mathbb{C} \text{ (or } \mathbb{R})$.

\medskip

Similarly, we introduce the set of places on function fields.

Let $C$ be a a smooth projective curve defined over an algebraically closed field $k$ and $L:=k(C)$ be the function field of $C$. The set of places on $L$, denoted by $\mathcal{M}_{L}$ consists of all absolute values of the form:  $$|\cdot|_p: x\mapsto e^{\ord_p(x)}$$ for any $x\in L$ and any closed point $p\in C$.

Let $K/L$ be a finite field extension. The set of places on $K$ is denoted by $\mathcal{M}_{K}$ and consists of all absolute values on $K$ whose restriction to $L$ is one of the places in $\mathcal{M}_{L}.$ In this case, all the places in $\mathcal{M}_K$ are nonarchimedean. Set $\mathcal{M}_K^0=\mathcal{M}_K$ and $\mathcal{M}_K^{\infty}=\emptyset$ for convenience.
\bigskip

Let $K/L$ be a finite field extension where $L=\mathbb{Q}$ or a function field $k(C)$ of a curve $C$. For any place $v\in \mathcal{M}_K$, denote by $n_v:=[K_v:L_v]$ the local degree of $v$ then we have the product formula $$\prod_{v\in \mathcal{M}_K}|x|_v^{n_v}=1$$  for all $x\in K^*$.

For any $v\in \mathcal{M}_K$, denote by $O_v:=\{x\in K|\,\, |x|_v\leq 1 \}$ the ring of $v$-integers. In the number field case, we also denote by $O_K:=\{x\in K|\,\, |x|_v\leq 1 \text{ for all }v\in \mathcal{M}_K^0\}$ the ring of integers.

\medskip

\subsection{Basics on Heights}We recall some basic properties of heights that are needed in the proof of Theorem A, see \cite{se1} or \cite{Silverman1986} for detail.

In this section, we set $L=\mathbb{Q}$ or $k(C)$ the function field of a curve $C$ defined over an algebraically closed field $k.$ Denote by $\overline{L}$ its algebraic closure.

\smallskip

\begin{prodefi}Let $K/L$ be a finite field extension. Let $p\in \mathbb{P}^n(K)$ be a point with homogeneous coordinate $p=[x_0:\cdots:x_n] \text{ where    }x_0,\cdots,x_n \in K.$ The $height$ of $p$ is the quantity $$H_{\mathbb{P}^n}(p):=(\prod_{v\in \mathcal{M}_K}\max\{|x_0|_v,\cdots,|x_n|_v \}^{n_v})^{1/[K:L]}.$$ The height $H_{\mathbb{P}^n}(p)$ depends neither on the choice of homogeneous coordinates of $p$, nor on the choice of a field extension $K$ which contains $p.$
\end{prodefi}

When $L=k(C)$, we have a geometric interpretation of the height $H_{\mathbb{P}^n}(p)$.
Observe that $\mathbb{P}^n_L$ is the generic fiber of the trivial fibration $\pi: \mathbb{P}^n_C:=\mathbb{P}^n\times C\rightarrow C.$ We set $s_p:D\rightarrow \mathbb{P}^n_C$ the normalization of the Zariski closure of $p$ in $\mathbb{P}^n_C.$
Then we have $$H_{\mathbb{P}^n}(p)=e^{\deg(s_p^{*}O_{\mathbb{P}^n_C}(1))/\deg(\pi\circ s_p)}.$$

\begin{pro}\label{prodh}Let $f:\mathbb{P}_{\overline{L}}^n\dashrightarrow \mathbb{P}_{\overline{L}}^m$ be a rational map and $X$ be a subvariety of $\mathbb{P}_{\overline{L}}^n$ such that $I(f)\bigcap X$ is empty and the restriction $f|_X$ is finite of degree $d$ onto its image $f(X)$.

Then there exist $A>0$ such that for all point $p\in X(\overline{L})$, we have $$\frac{1}{A}H_{\mathbb{P}^n}(p)^d\leq H_{\mathbb{P}^m}(f(p))\leq AH_{\mathbb{P}^n}(p)^d.$$
\end{pro}

\begin{pro}[Northcott Property]Let $K/\mathbb{Q}$ be a number field, and $B>0$ be any constant. Then the set $$\{p\in \mathbb{P}^n(K)|\,\, H_{\mathbb{P}^n}(p)\leq B\}$$ is finite.
\end{pro}
\rem The Northcott Property does not hold in the case $K=k(C)$ when $k$ is not a finite field. For example, the set $$\{p\in \mathbb{P}^n(k(t))|\,\, H_{\mathbb{P}^n}(p)=0\}=\{[x:y]|\,\,(x,y)\in k^2\setminus \{(0,0)\}\}$$ is infinite.
\endrem

\subsection{Upper bounds on heights}\label{uboh}

Let $K$ be a number field or a function field of a smooth curve over an algebraically closed field $k'$.
Let $f:\mathbb{A}^2_{\overline{K}}\rightarrow\mathbb{A}^2_{\overline{K}}$ be any birational polynomial morphism defined over $\overline{K}$ and assume that $\lambda(f)>1$.

 According to Theorem \ref{l>1goodmodel}, we may suppose that there exists a compactification $X$ of $\mathbb{A}^2_{\overline{K}}$
, a closed point $Q\in X\setminus \mathbb{A}^2_{\overline{K}}$ such that $f$ extends to a birational transformation $\widetilde{f}$ on $X$ which satisfies the following properties:
\begin{points}
\item
$\widetilde{f}$ is algebraically stable on $X$;
\item
there exists a closed point $Q\in X \setminus \mathbb{A}^2$ fixed by $\widetilde{f}$, such that $d\widetilde{f}(Q)=0$;
\item
$\widetilde{f}(X \setminus \mathbb{A}^2)=Q.$
\end{points}
To simplify, we write $f=\widetilde{f}$ in the rest of the paper.
We fix an embedding $X\subseteq \mathbb{P}^N_{\overline{K}}.$ Let $C$ be an irreducible curve in $X$ whose intersection with $\mathbb{A}^2_{\overline{K}}$ is non empty.

\begin{pro}\label{lemupbound} Suppose that $C$ is not periodic and $C\setminus \mathbb{A}^2_{\overline{K}}=\{Q\}.$
 Then there exists a number $B>0$ such that for any point $p\in C(K)$ for which
the set $\{n\in \mathbb{N}|\,\, f^n(p)\in C\}$ is infinite, we have $H_{\mathbb{P}^N}(p)\leq B.$
\end{pro}
\proof
Assume that $X, f, C$ and $Q$ are all defined over $K$ and $Q=[1:0:\cdots: 0]\in \mathbb{P}^N_{K}$. We can extend $f$ to a rational morphism on $\mathbb{P}^N$ which is regular at $Q$. Then there exists an element $a\in K^*$ and $F_i\in (x_1,\cdots,x_N)K[x_0,\cdots,x_N]$ for $i=0,\cdots , N$ such that $$f([1:x_1:\cdots:x_N])=[a+F_0:F_1:\cdots :F_N]$$ for any $[1:x_1:\cdots: x_N]\in X.$ Since $f$ is regular at $Q$ and $a\neq 0,$ there is a finite set $S\subseteq \mathcal{M}_K^0$ such that for any $v\in \mathcal{M}_K^0\setminus S$, we have $|a|_v=1$ and all coefficients of $f$ are defined in $O_v$. Recall that we may endow $X$ with a metric $d_v$, see Section \ref{secmop}.

For any $v\in \mathcal{M}_K^0\setminus S$, set $r_v:=1$ and $U_v:=\{x\in X(K)|\,\, d_v(x,Q)<1\}.$ Since $df(Q)=0$, we see that for all $x\in U_v$, $d_v(f(x),Q)\leq d_v(x,Q)^2$, hence $$\lim_{n\rightarrow \infty}f^n(x)=Q.$$

For any $v\in S$, set $r_v:=|a|_v$ and $U_v:=\{x\in X(K)|\,\, d_v(x,q)<r_v\}.$ We see that for all $x\in U_v$, $d_v(f(x),Q)\leq d_v(x,Q)^2/r_v$, and again it follows that $$\lim_{n\rightarrow \infty}f^n(x)=Q.$$

For any $v\in \mathcal{M}_K^{\infty}$, since $df(Q)=0$, there is $r_v>0$ such that for any $x\in U_v:=\{x\in X(K)|\,\, d_v(x,q)<r_v\}$ we have $f(x)\subseteq U_v$ and $$\lim_{n\rightarrow \infty}f^n(x)=Q.$$

If $p\in \bigcup_{v\in \mathcal{M}_K}U_v$, Theorem \ref{local} shows that $C$ is periodic and this contradicts our assumption.
In other words, we need to estimate the height of a given point $$p\in C(K)\setminus \bigcup_{v\in \mathcal{M}_K}U_v.$$
If $C$ intersects the line at infinity only at the point $Q$, then we may directly
estimate the height of $p$ given by the embedding of $C$ into $\mathbb{P}^N_K$. Since we do not
assume that this is the case, we shall work first with a height induced by a divisor on $C$ given by the divisor $Q$,
and then estimate
$h_{\mathbb{P}^N}( p )$ using  Proposition  \ref{prodh}.
To do so, let $i:\widetilde{C}\rightarrow C\subseteq X$ be the normalization of $C$ and pick a point $Q'\in i^{-1}(Q)$. There is a positive integer $l$ such that $lQ'$ is a very ample divisor of $\widetilde{C}$. So there is an embedding $j: \widetilde{C}\hookrightarrow \mathbb{P}^M$ for some $M>0$ such that $$Q'=[1:0:\cdots :0]=H_{\infty}\bigcap \widetilde{C}$$  where $H_{\infty}=\{x_M=0\}$ is the hyperplane at infinity. Let $$g:\widetilde{C}\rightarrow \mathbb{P}^1$$ be a morphism sending $[x_0:\cdots :x_M]\in \widetilde{C}$ to $[x_0:x_M]\in \mathbb{P}^1$. It is well defined since $\{x_0=0\}\bigcap H_{\infty}\bigcap \widetilde{C}=\emptyset$. Then $g$ is finite and $$g^{-1}([1:0])=H_{\infty}\bigcap\widetilde{C}=[1:0\cdots:0].$$
By base change, we may assume that $\widetilde{C},i,j, g$ are all defined over $K.$

In the function field case, there is a smooth projective curve $D$ such that $K=k'(D)$; and in the number field case, we set $D=\Spec O_K.$

 We consider the irreducible scheme $\widetilde{\mathcal{C}}\subseteq \mathbb{P}^M_{D}$ over $D$ whose generic fiber is $\widetilde{C}$ and the irreducible scheme $\mathcal{X}\subseteq \mathbb{P}^N_{D}$ over $D$ whose generic fiber is $X.$
 Then $i$ extends to a map $\iota:\widetilde{\mathcal{C}}\dashrightarrow \mathcal{X}$ over $D$ birationally to its image.
For any $v\in \mathcal{M}^0_K$, let $$\mathfrak{p}_v=\{x\in O_v|\,\,v(x)>0\}$$ be a prime ideal in $O_v.$ There is a finite set $T$ consisting of those places $v\in \mathcal{M}^0_K$ such that $\iota$ is not regular along the special fibre $\widetilde{C}_{O_v/\mathfrak{p}_v}$ at $\mathfrak{p}_v\in D$ or $\widetilde{C}_{O_v/\mathfrak{p}_k}\bigcap H_{\infty, O_v/\mathfrak{p}_v}\neq \{[1:0:\cdots:0]\}$.

For any $v\in \mathcal{M}^0_K-T\bigcup S,$ observe that we have
 \begin{align*}
V_v:=&\{[1:x_1:\cdots:x_M]\in \widetilde{C}(K)|\,\, |x_i|_v<1 ,i=1,\cdots, M\}\nonumber\\
         =&\{[1:x_1:\cdots:x_M]\in \widetilde{C}(K)|\,\, |x_M|_v<1 \}=g^{-1}(\Omega_v)\bigcap \widetilde{C}(K)
\end{align*}
with $\Omega_v:=\{[1:x]\in \mathbb{P}^1(K)|\,\,|x|_v<t_v\}$ and $t_v:=1$.

For any $v\in T\bigcup S\bigcup \mathcal{M}^{\infty}_K$, by the continuity of $i$, there is $s_v>0$ such that $$i(V_v)\in U_v$$ where $V_v=\{[1:x_1:\cdots:x_M]\in \widetilde{C}(K)|\,\,|x_i|_v<s_v ,i=1,\cdots M\}.$  Since $g^{-1}([1:0])=\{[1:0:\cdots :0]\},$ there exists $t_v>0$, such that $$g^{-1}(\Omega_v)\bigcap \widetilde{C}\subseteq V_v$$ where $\Omega_v=\{[1:x]\in \mathbb{P}^1(K)|\,\,|x|_v<t_v\}.$

\medskip

We need to find an upper bound for the height of points in $C(K)\setminus\bigcup_{v\in \mathcal{M}_K}U_v.$ Since the set $\Sing(C)$ of singular points of $C$ is finite, we only have to bound the height of points in $C(K)\setminus(\Sing(C)\bigcup_{v\in \mathcal{M}_K}U_v).$

Let $p$ be a point in $C(K)\setminus(\Sing(C)\bigcup_{v\in \mathcal{M}_K}U_v)$. Observe that $i^{-1}(p)\in \widetilde{C}(K)$ and $x:=j(i^{-1}(p))$ is also defined over $K$.  We have $x\not\in V_v$ hence $y:=g(x)\not\in \Omega_v$ for all $v\in \mathcal{M}_K.$

For any $y=[y_0:y_1]\in \mathbb{P}^1(K)\setminus(\bigcup_{v\in \mathcal{M}_K}\Omega_v)$, we have $|y_1/y_0|_v\geq t_v$ for all $v.$ We get the following upper bound
\begin{align*}
H_{\mathbb{P}^1}(y)^{[K:\mathbb{Q}]}=&\prod_{v\in \mathcal{M}_K}\max\{|y_0|_v,|y_1|_v\}^{n_v}\nonumber\\
         \leq&\prod_{v\in\mathcal{M}_K}\max\{|y_1|_v/t_v,|y_1|_v\}^{n_v}\\
            =&\prod_{v\in \mathcal{M}_K}|y_1|_v^{n_v}\prod_{v\in \mathcal{M}_K}\max\{1,1/t_v\}^{n_v}\\
            =&\prod_{v\in \mathcal{M}_K}\max\{1,1/t_v\}^{n_v}=: B'<\infty.
\end{align*}
By Proposition \ref{prodh} applied to $g:\widetilde{C}\hookrightarrow \mathbb{P}^M$ and $i:\widetilde{C}\rightarrow \mathbb{P}^N$, we get $H_{\mathbb{P}^N}(p)\leq B$ for some constant $B$ independent on the choice of $p$ as we require.
\endproof

\section{Proof of Theorem A}\label{seproofa}
 Let $C$ be a curve in $\mathbb{A}^2_{K}$.
We want to show that  for any point $p\in \mathbb{A}^2(K)$ such that the set $$\{n\in \mathbb{N}|\,\, f^n(p)\in C\}$$ is infinite, then $p$ is preperiodic.

 According to Section \ref{l=1}, we may assume that $\lambda(f)>1.$
 As in Section \ref{uboh}, we use
Theorem \ref{l>1goodmodel} to get a compactification $X$ of $\mathbb{A}^2_{K}$. For simplicity, we also denote by $f$ the map induced by $f$ on $X.$ There exists $n\geq 1$ such that $f^n$ contracts $X\setminus \mathbb{A}_{K}^2$ to a superattracting fixed point $Q\in X\setminus\mathbb{A}_{K}^2.$ We extend $C$ to a curve in $X$. Suppose that $C$ is not periodic. By Theorem \ref{qinf^nC}, we may assume that $C(K)\setminus\mathbb{A}^2(K)=\{Q\}.$ Finally we
fix an embedding $X\hookrightarrow \mathbb{P}^N_{K}$ for some $N\geq 1.$
\bigskip

We first treat the case $K=\overline{\mathbb{Q}}$.

There is a number field $K'$ such that both $f$ and $p$ are defined over $K'.$ Then $f^n(p)\in \mathbb{A}^2(K')$ for all $n\geq 0.$

Proposition \ref{lemupbound} and the Northcott Property imply that the set $\{f^n(p)|\,\, n\geq 0\}\bigcap C$ is finite. Since the set $\{n\in \mathbb{N}| \,\, f^n(p)\in C\}$ is infinite, there exists $n_1>n_2>0$ such that $f^{n_1}(p)=f^{n_2}(p).$ We conclude that $p$ is preperiodic.
\medskip

Next we consider the general case of an algebraically closed field $K$ of characteristic $0.$

By replacing $K$ by an algebraically closed subfield over which $p,C$ and $f$ are all defined, we may suppose that the transcendence degree $\trd K/\mathbb{Q}$ of $K$ over $\mathbb{Q}$ is finite.
 We argue by induction on $\trd K/\mathbb{Q}$.

If $\trd K/\mathbb{Q}=0$, then $K=\overline{\mathbb{Q}}$ and we are due by what precedes.

If $\trd K/\mathbb{Q}\geq 1$, then there is an algebraically closed subfield $k$ of $K$ such that $\trd k/\mathbb{Q}=\trd K/\mathbb{Q}-1.$

There is a  smooth projective curve $D$ over $k$, such that $X,f,p, Q$ and $C$ are defined over the function field $k(D)$ of $D$. Observe that $K=\overline{k(D)}.$

 We consider the irreducible scheme $$\pi:\mathcal{X}\subseteq \mathbb{P}^N_{D}\rightarrow D$$ over $D$ whose generic fiber is $X$ and $\mathcal{C}\subseteq \mathbb{P}^N_{D}$ the Zariski closure of $C$ in $\mathcal{X}.$

 The map $f$ extends to a birational map $f':\mathcal{X}\dashrightarrow \mathcal{X}$ over $D$. For any $x\in D$, denote by $X_x$ and $C_x$ the fiber of $\mathcal{X}$ and $\mathcal{C}$ at $x$ respectively, and denote by $f_x$ the restriction of map $f'$ to the fiber $X_x$.

Proposition \ref{lemupbound} implies that there is a number $M\geq 0$ such that for all $n\geq 0$ either $f^n(p)\not\in C$ or
$H_{\mathbb{P}^N}(f^n(p))\leq M.$

A point $s\in X(k(D))$ is associated to its Zariski closure in $\mathcal{X}$ which is a section of $\pi:\mathcal{X}\rightarrow D$. For simplicity, we also write $s$ for this section. Then the height of $s$ is $$H_{\mathbb{P}^N}(s)=e^{(s\cdot L)}$$ where $L:=O_{\mathbb{P}_D^N}(1)$.

For any section $s$, observe that $\pi$ induces an isomorphism from $s$ to the curve $D$. We may consider the Hilbert polynomial $$\chi(L^{\otimes n},s)=1-g(s)+n(s\cdot L)=1-g(D)+n\log H(s).$$ It follows that there is a quasi-projective $k$-variety $M_H$ that parameterizes the sections $s$ of $\pi$ such that $H_{\mathbb{P}^N}(s)\leq M$ (see \cite{debarre}).

Let $T_1$ be the set of points  $x\in D$ such that $f_x$ is birational and $I(f_x^{-1})\bigcap I(f_x)\neq \emptyset$. Observe that $T_1$ is finite.  Let $T_2$ be the set of the points $x\in D\setminus T_1$, such that $C_x$ is fixed. Since $C$ is not fixed, $T_2$ is finite. Because $k$ is algebraically closed, $D\setminus(T_1\bigcup T_2)$ is infinite.  For any point $x\in D$, denote by $p_x:M_H\rightarrow X_x$ the map sending $s$ to $s(x)$.  Pick a sequence of distinct points $\{x_i\}_{i\geq 0}\subseteq D\setminus(T_1\bigcup T_2)$. For any $l\geq 1$, let $$p_l=\prod_{i=1}^{l}p_{x_i}: M_H\rightarrow \prod_{i=1}^{l}X_{x_i}.$$ Observe that any two points $s_1,s_2\in M_H$ are equal if and only if $p_i(s_1)= p_i(s_2)$ for all $i\geq 0$.

We claim the following lemma, and prove it later.
\begin{lem}\label{embeding}Let $X$ be any reduced quasi-projective variety over an algebraically closed field $k$. For any $i\geq 1$, let $\pi_i: X\rightarrow Y_i$ be a morphism. If for any difference points $x_1,x_2\in X$, there exists $i\geq 0$, such that $\pi_i(x_1)\neq \pi_i(x_2)$, then for $l$ large enough the map $$p_l=\prod_{i=1}^{l}\pi_{i}: X\rightarrow \prod_{i=1}^{l}Y_i$$ is finite.
\end{lem}
By Lemma \ref{embeding}, there is an integer $L$ large enough, such that the map $p_L$ is finite. By Proposition \ref{periodictofixed}, $C_{x_i}$ is not periodic for all $i\geq 1.$ The set $\mathsf{N}:=\{n\geq 0|f^n(p)\in C\}$ is infinite, enumerate $\mathsf{N}=\{n_1<n_2<\cdots<n_i<n_{i+1}<\cdots\}$. For any $i\geq 0,$ there exists $s_i\in M_H$ such that $s_i=f^{n_i}(p).$  By the induction hypothesis, we know that $s_{n_0}(x_i)=f^{n_0}(p)(x_i)$ is a preperiodic point of $f_{x_i}$ for any $1\leq i \leq L$. Then the orbit $G_i$ of $p(x_i)$ in $X_{x_i}$ is finite. So the set $$p_L(\{s_i\}_{i\geq 0})\subseteq \prod_{i=0}^L G_i$$ is finite. Since $p_L$ is finite, then we have $\{s_i\}_{i\geq 0}$ is finite. Then there is $i_1>i_2$ such that $s_{i_1}=s_{i_2}$, and $f^{n_{i_1}}(p)=f^{n_{i_2}}(p)$. Then $p$ is preperiodic.
\endproof
\proof[Proof of Lemma \ref{embeding}] We prove this lemma by induction on the dimension of $X.$

If $\dim X=0$, then the result is trivial.

If $\dim X>0$, we may assume that $X$ is irreducible. We pick any point $x\in X$, and let $F_l$ be the fiber of $p_l$ which contains $x.$ Observe that $$F_{l+1}\subseteq F_{l},$$ so that there is an integer $L'\geq 1$, such that for any $L\geq L'$, $$F_L=\bigcap_{l\geq 0}F_l.$$ Since for any point $x_1\in X-\{x\}$, there exists $i\geq 0$, such that $\pi_i(x_1)\neq \pi_i(x)$, we have $$F_L=\bigcap_{l\geq 0}F_l=\{x\},$$ so that $$\dim X-\dim p_L(X)\leq \dim F_L=0.$$ In particular $p_L$ is generically finite. It means that there exists an open set $U$ of $p_L(X)$, such that $p_L: p_L^{-1}(U)\rightarrow U$ is finite. Set $X'=X-p_L^{-1}(U)$, then we have $\dim X'\leq \dim X-1.$

By the induction hypothesis, there is $L''\geq L'$, such that for any $L\geq L''$, $p_L|_{X'}$ is finite and then $p_L$ is finite.
\endproof

\bibliography{dd}

\end{document}